\documentclass[11pt]{amsart}
\usepackage{eucal,fullpage,times,amsmath,amsthm,amssymb,mathrsfs,stmaryrd}
\usepackage[all]{xy}
\usepackage{url}

\newcommand{\calO}{{\mathcal{O}}}
\newcommand{\calL}{{\mathcal{L}}}
\newcommand{\calU}{{\mathcal{U}}}

\newcommand{\calA}{{\mathcal{A}}}

\newcommand{\A}{\mathbf{A}}
\newcommand{\G}{\mathbf{G}}
\newcommand{\Z}{\mathbf{Z}}
\newcommand{\C}{\mathbf{C}}
\newcommand{\F}{\mathbf{F}}
\newcommand{\Q}{\mathbf{Q}}
\renewcommand{\P}{\mathbf{P}}

\newcommand{\Spec}{{\mathrm{Spec}}}

\newcommand{\Hom}{\mathrm{Hom}}
\newcommand{\Map}{\mathrm{Map}}

\newcommand{\Sym}{\mathrm{Sym}}

\newcommand{\R}{\mathrm{R}}

\newcommand{\et}{\mathrm{\acute{e}t}}
\newcommand{\fppf}{\mathrm{fppf}}

\newcommand{\zar}{\mathrm{Zar}}
\newcommand{\prop}{\mathrm{prop}}
\newcommand{\im}{\mathrm{im}}

\newcommand{\fm}{\mathfrak{m}}
\newcommand{\fn}{\mathfrak{n}}
\newcommand{\fp}{\mathfrak{p}}

\newcommand{\comment}[1]{}

\begin{document}

\bibliographystyle{alpha}

\newtheorem{theorem}{Theorem}[section]
\newtheorem*{theorem*}{Theorem}
\newtheorem*{condition*}{Condition}
\newtheorem*{definition*}{Definition}
\newtheorem{proposition}[theorem]{Proposition}
\newtheorem{lemma}[theorem]{Lemma}
\newtheorem{corollary}[theorem]{Corollary}
\newtheorem{claim}[theorem]{Claim}
\newtheorem{claimex}{Claim}[theorem]

\theoremstyle{definition}
\newtheorem{definition}[theorem]{Definition}
\newtheorem{question}[theorem]{Question}
\newtheorem{remark}[theorem]{Remark}
\newtheorem{example}[theorem]{Example}
\newtheorem{condition}[theorem]{Condition}
\newtheorem{warning}[theorem]{Warning}
\newtheorem{notation}[theorem]{Notation}

\title{Annihilating the cohomology of group schemes}
\author{Bhargav Bhatt}
\begin{abstract} Our goal in this note is to show that cohomology classes with coefficients in finite flat group schemes can be killed by finite covers of the base scheme, and similarly for abelian schemes with ``finite covers'' replaced by ``proper covers.'' We apply this result to commutative algebra by giving a new and more conceptual proof of Hochster-Huneke's theorem on the existence of big Cohen-Macaulay algebras in positive characteristic; all previous proofs of this result were equational or cocycle-theoretic in nature.
\end{abstract}

\maketitle

\section{Introduction}

Given a scheme $S$ with a sheaf $G$ and class $\alpha \in H^n(S,G)$ for $n > 0$, a natural question one may ask is if there exist covers $\pi:T \to S$ such that $\pi^* \alpha = 0$? Of course, as stated, the answer is trivially yes as we may take $T$ to be a disjoint union of suitable opens occurring in a Cech cocyle representing $\alpha$. However, the question becomes interesting if we require geometric conditions on $\pi$, such as properness or even finiteness. Our goal is to study such questions for fppf cohomology in the case that $G$ is either a finite flat commutative group scheme or an abelian scheme.  Our main results are:

\begin{theorem}
\label{killcohffgs}
Let $S$ be a noetherian excellent scheme, and let $G$ be a finite flat commutative group scheme over $S$. Then classes in $H^n_\fppf(S,G)$ can be killed by finite surjective maps to $S$ for $n > 0$.
\end{theorem}

\begin{theorem}
\label{killcohabs}
Let $S$ be a noetherian excellent scheme, and let $A$ be an abelian scheme over $S$. Then classes in $H^n_\fppf(S,A)$ can be killed by proper surjective maps to $S$ for $n > 0$. Moreover, there exists an example of a normal affine scheme $S$ that is essentially of finite type over $\C$, and an abelian scheme $A \to S$ with a class in $H^1_\fppf(S,A)$ that cannot be killed by finite surjective maps to $S$.
\end{theorem}

Our primary motivation for proving the preceding results was to obtain a better understanding of the Hochster-Huneke proof of the existence of big Cohen-Macaulay algebras in positive characteristic commutative algebra (see \cite{HHBigCM}). We have succeeded in this endeavour as we can give a new and essentially {\em topological} proof of the Hochster-Huneke result by using the cohomology-annihilation results discussed above in lieu of the more tradtional equational approaches; see \S \ref{sec:altproof} for more. We are hopeful that a similar approach, coupled with Fontaine's constructions of mixed characteristic rings admitting Frobenius actions (see \cite{Fontainepadicperiod}), will eventually provide an approach to Hochster's homological conjectures in mixed characteristic commutative algebra; we refer the interested reader to \cite{HochsterSurvey} for further information.

\subsubsection*{An informal summary of the proofs:} To prove Theorem \ref{killcohffgs}, we first use a theorem of Raynaud to embed a finite flat group scheme into an abelian scheme; this permits a reduction to from fppf cohomology to \'etale cohomology by a theorem of Grothendieck. Next, using an observation due to Gabber, we reduce from \'etale cohomology to Zariski cohomology, and then we solve the problem by hand. For Theorem \ref{killcohabs}, we reduce as before to Zariski cohomology, and then solve the problem using de Jong's alterations results combined with an observation concerning rational sections of an abelian scheme over a regular base scheme. The example referred to in Theorem \ref{killcohabs} is constructed using the cone on an elliptic curve and a geometric idea; the existence of such an example is suggested by analysing the proof of the first half of Theorem \ref{killcohabs}. Lastly, the Hochster-Huneke theorem is reproven by by first reformulating it as a suitable cohomology-annihilation statement for the higher local cohomology of the structure sheaf, and then deducing this statement from Theorem \ref{killcohffgs} by using finite flat subgroup schemes of $\G_a$ defined by additive polynomials in Frobenius.
\subsubsection*{Notations and conventions}
All group schemes occurring in this note are commutative; all the cohomology groups occurring in this note are computed in the fppf topology unless otherwise specified. We stress that there are no assumptions on the residue characteristics of the base scheme $S$ in Theorems \ref{killcohffgs} and \ref{killcohabs}. 

\subsubsection*{Organisation of this note}
In \S \ref{sec:gabberresult} we recall Gabber's observation alluded to above. Using this observation, we prove Theorem \ref{killcohffgs} in \S \ref{sec:ffgsresult}, and the first half of Theorem \ref{killcohabs} in \S \ref{sec:abschresult}. Next, in \S \ref{sec:altproof}, we explain how to use Theorem \ref{killcohffgs} to give a new proof of the Hochster-Huneke theorem. We close in \S \ref{abspurityneedh} by giving an example that illustrates the necessity of ``proper'' in the first half of Theorem \ref{killcohabs} and finishes its proof. 

\subsubsection*{Acknowledgements:} The results presented here formed a part of my doctoral thesis written under Johan de Jong, and I would like to thank him for stimulating discussions. I would also like to thank J\'anos Koll\'ar for a very useful conversation.

\section{An observation of Gabber}
\label{sec:gabberresult}

In this section, we recall a result of Gabber concerning the local structure of the \'etale topology. This observation permits reduction of \'etale cohomological considerations to those in finite flat cohomology and those in Zariski cohomology. We begin with an elementary lemma on extending covers that will be used repeatedly in the sequel.

\begin{lemma}
\label{extendcover}
Fix a noetherian scheme $X$. Given an open dense subscheme $U \to X$ and a finite (surjective) morphism $f:V \to U$, there exists a finite (surjective) morphism $\overline{f}:\overline{V} \to X$ such that $\overline{f}_U$ is isomorphic to $f$. Given a Zariski open cover $\calU = \{j_i:U_i \to X\}$ with a finite index set, and finite (surjective) morphisms $f_i:V_i \to U_i$, there exists a finite (surjective) morphism $f:Z \to X$ such that $f_{U_i}$ factors through $f_i$. The same claims hold if ``finite (surjective)'' is replaced by ``proper (surjective)'' everywhere.
\end{lemma}
\begin{proof}
We first explain how to deal with the claims for finite morphisms. For the first part, Zariski's main theorem (\cite[Th\'eor\`eme 8.12.6]{EGA4_3}) applied to the morphism $V \to X$ gives a factorisation $V \hookrightarrow W \to X$ where $V \hookrightarrow W$ is an open immersion, and $W \to X$ is a finite morphism. The scheme-theoretic closure $\overline{V}$ of $V$ in $W$ provides the required compactification in view of the fact that finite morphisms are closed.

For the second part, by the above, we may extend each $j_i \circ f_i:V_i \to X$ to a finite surjective morphism $\overline{f_i}:\overline{V_i} \to X$ such that $\overline{f_i}$ restricts to $f_i$ over $U_i \hookrightarrow X$. Setting $W$ to be the fibre product over $X$ of all the $\overline{V_i}$ is then seen to solve the problem.

To deal with the case of proper (surjective) morphisms instead of finite (surjective) , we repeat the same argument as above replacing the reference to Zariski's main theorem by one to Nagata's compactification theorem (see \cite[Theorem 4.1]{ConradNagata}).
\end{proof}

Next, we sketch Gabber's result.

\begin{lemma}[{Gabber, \cite[Lemma 5]{HooblerGabber}}]
\label{gabberetcovaff}
Let $f:U \to X$ be a surjective \'etale morphism of affine schemes. Then there exists a finite flat map $g:X' \to X$, and a Zariski open cover $\{U_i \hookrightarrow X'\}$ such that the natural map $\sqcup_i U_i \to X$ factors through $U \to X$.
\end{lemma}

For the convenience of the reader, we sketch a proof.

\begin{proof}[Sketch of proof]
We first explain how to deal with the local case. Assume that $X = \Spec(A)$ is the spectrum of a local ring $A$, and $U = \Spec(B)$ is the spectrum of a local \'etale $A$-algebra $B$. The structure theorem for \'etale morphisms (see \cite[Expos\'e I, Th\'eor\`eme 7.6]{SGA1}) implies that $B = C_\fm$ where $C = A[x]/(f(x))$ with $f(x) = x^n + a_1 x^{n-1} + \cdots + a_n$ a monic polynomial, and $\fm \subset C$ a maximal ideal with $f'(x) \notin \fm$. We define 
\[D = A[x_1,\cdots,x_n]/(\sigma_i(x_1,\cdots,x_n) - (-1)^{n-i} a_i) \]
where $\sigma_1,\cdots,\sigma_n$ are the elementary symmetric polynomials in the $x_i$'s. This ring is finite free over $A$ of rank $n!$, admits an action of $S_n$ that is transitive on the maximal ideals, and formalises the idea that the coefficients of $f(x)$ can be written as elementary symmetric functions in its roots. In particular, there is a natural morphism $C \to D$ sending $x$ to $x_1$. As both $C$ and $D$ are finite free over $A$, there is a maximal ideal $\fm_1 \subset D$ lying over $\fm \subset C$. Thus, there is a natural map $a:B \to D_{\fm_1}$. By the $S_n$-action, for every maximal ideal $\fn \subset D$, there is an automorphism $D \to D$ sending $\fm_1$ to $\fn$. Composing such an automorphism with $a$, we see that for every maximal ideal $\fn \subset D$, the structure map $A \to D_\fn$ factorises through $A \to B$ for some map $B \to D_\fn$; the claim follows. In general, one reduces to the local case by considering fibre products as in Lemma \ref{extendcover}
\end{proof}

Actually, we use a slight weakening of Gabber's result -- relaxing finite flat to finite surjective -- that remains true when the schemes under consideration are no longer assumed to be affine.

\begin{lemma}
\label{gabberetcov}
Let $f:U \to X$ be a surjective \'etale morphism of schemes. Then there exists a finite surjective map $g:X' \to X$, and a Zariski open cover $\{U_i \hookrightarrow X'\}$ such that the natural map $\sqcup_i U_i \to X$ factors through $U \to X$.
\end{lemma}
\begin{proof}
We can solve the problem locally on $X$ by Lemma \ref{gabberetcovaff}. This means that there exists a Zariski open cover $\{V_i \hookrightarrow X\}$, finite surjective maps $W_i \to V_i$, and Zariski covers $\{Y_{ij} \hookrightarrow W_i\}$ such that $\sqcup Y_{ij} \to V_i$ factors throughts $U \times_X V_i \to V_i$. By Lemma \ref{extendcover}, we may find a single finite surjective map $W \to X$ such that $W \times_X V_i \to V_i$ factors through $W_i \to V_i$. Setting $X' = W$ and  pulling back the covers $\{Y_{ij} \to W_i\}$ to $W \times_X V_i$ then solves the problem.
\end{proof}

\section{The theorem for finite flat commutative group schemes}
\label{sec:ffgsresult}

In this section we prove Theorem \ref{killcohffgs} following the plan explained in the introduction. To carry that program out, we first explain how to relate the fppf cohomology of finite flat group schemes to \'etale cohomology; it turns out that they are almost the same.

\begin{proposition}
\label{ffgscohshring}
Let $S$ be the spectrum of a strictly henselian local ring, and let $G$ be a finite flat commutative group scheme over $S$. Then $H^i(S,G) = 0$ for $i > 1$.
\end{proposition}
\begin{proof}
We first explain the idea informally. Using a theorem of Raynaud, we can embed $G$ into an abelian scheme, which allows us to express the cohomology of $G$ in terms of that of abelian schemes. As abelian schemes are smooth, a result of Grothendieck ensures that their fppf cohomology coincides with their \'etale cohomology. As the latter vanishes when $S$ is strictly henselian, we obtain the desired conclusion.

Now for the details: a construction of Raynaud (see \cite[Th\'eor\`eme 3.1.1]{BBMII}) gives the existence of an abelian scheme $A \to S$ and an $S$-closed immersion $G \hookrightarrow A$ of group schemes. Let $A/G$ denote the quotient sheaf for the fppf topology. It is well known that $A/G$ is an abelian scheme over $S$, but we give a proof here for the lack of a suitable reference; the reader is advised to skip to the next paragraph. A theorem of Michael Artin (see \cite[Corollary 6.3]{ArtinVersalDefs}) ensures that $A/G$ is an algebraic stack, i.e., admits a smooth presentation. Since the action of $G$ on $A$ is free, the quotient is a sheaf of discrete groupoids and hence actually an algebraic space (see \cite[Corollary 10.4]{LMBChAlg}). Since $A$ and $G$ are locally of finite presentation implies the same is true for $A/G$. Moreover, the quotient map $A \to A/G$ is flat since it is a $G$-torsor. The Auslander-Buschbaum theorem then forces the geometric fibres of $A/G \to S$ to be regular. The fibre-by-fibre smoothness criterion (see \cite[Th\'eor\`eme 17.5.1]{EGA4_4}) shows that $A/G \to S$ is smooth.  The map $A/G \to S$  is proper with geometrically connected fibres by descent since the same is true for $A \to S$. Hence, we find that $A/G$ is an algebraic space over $S$ representing a sheaf of abelian groups with the additional property that the structure morphism $A/G \to S$ is a proper smooth morphism with geometrically connected fibre. A theorem of Raynaud (see \cite[Theorem 1.9]{FaltingsChai}), then shows that $A/G$ is actually an abelian scheme over $S$. 

The preceding construction gives us a short exact sequence
\[ 0 \to G \to A \to A/G \to 0 \]
of abelian sheaves on the fppf site of $S$ relating the finite flat commutative group scheme $G$ to the abelian schemes $A$ and $A/G$. This gives rise to a long exact sequence
\[ \cdots H^{n-1}(S,A/G) \to H^n(S,G) \to H^n(S,A) \to \cdots \]
of fppf cohomology groups. By Grothendieck's theorem (see \cite[Th\'eor\`eme 11.7]{GrothDixExpIII}), fppf cohomology coincides with \'etale cohomology when the coefficients are smooth group schemes. Applying this to $A$ and $A/G$ then shows $H^i(S,A) = H^i(S,A/G) = 0$ for $i > 0$ as $S$ is strictly henselian. The claim about $G$ now follows from the preceding exact sequence.
\end{proof}

Next, we explain how to deal with Zariski cohomology with coefficients in a finite flat group scheme.

\begin{proposition}
\label{killcohzarffgs}
Let $S$ be a normal noetherian scheme, and let $G \to S$ be a finite flat commutative group scheme. Then $H^n_\zar(S,G) = 0$ for $n > 0$.
\end{proposition}
\begin{proof}
We may assume that $S$ is connected. As constant sheaves on irreducible topological spaces are acyclic, it will suffice to show that $G$ restricts to a constant sheaf on the small Zariski site of $S$, i.e., that the restriction maps $G(S) \to G(U)$ are bijective for any non-empty open subset $U \hookrightarrow S$. Injectivity follows from the density of $U \hookrightarrow S$ and the separatedness of $G \to S$. To show surjectivity, we note that given a section $U \to G$ of $G$ over $U$, we can simply take the scheme-theoretic closure of $U$ in $G$ to obtain an integral  closed subscheme $S' \hookrightarrow G$ such that the projection map $S' \to S$ is finite and an isomorphism over $U$. By the normality of $S$, this forces $S' = S$. Thus, $G$ restricts to a constant sheaf on $S$, as claimed.
\end{proof}

We can now complete the proof of Theorem \ref{killcohffgs} by following the outline sketched in the introduction.

\begin{proof}[Proof of Theorem \ref{killcohffgs}]
Let $S$ be a noetherian excellent scheme, and let $G \to S$ be a finite flat commutative group scheme. We need to show that classes in $H^n(S,G)$ can be killed by finite covers for $S$ for $n > 0$. We deal with the $n = 1$ case on its own, and then proceed inductively.

For $n = 1$, note that classes in $H^1(S,G)$ are represented by fppf $G$-torsors $T$ over $S$. By faithfully flat descent for finite flat morphisms, such schemes $T \to S$ are also finite flat. Passing to the total space of $T$ trivialises the $G$-torsor $T$. Therefore, classes in $H^1(S,G)$ can be killed by finite flat covers of $S$.

We now fix an integer $n > 1$ and a cohomology class $\alpha \in H^n(S,G)$. By Proposition \ref{ffgscohshring}, we know that there exists an \'etale cover of $S$ over which $\alpha$ trivialises. By Lemma \ref{gabberetcov}, after replacing $S$ by a finite cover, we may assume that there exists a Zariski cover $\calU = \{U_i \hookrightarrow S\}$ such that $\alpha|_{U_i}$ is Zariski locally trivial. The Cech spectral sequence for this cover is
\[ H^p(\calU,H^q(G)) \Rightarrow H^{p+q}(S,G) \]
where $H^q(G)$ is the Zariski presheaf $V \mapsto H^q(V,G)$. By construction, the class $\alpha$ comes from some $\alpha' \in H^{n-q}(\calU,H^q(G))$ with $q < n$.  The group $H^{n-q}(\calU,H^q(G))$ is the $(n-q)$-th cohomology group of the standard Cech complex
\[ \prod_i H^q(U_i,G) \to \prod_{i<j} H^q(U_{ij},G) \to \dots \]
By the inductive assumption and the fact that $q < n$, terms of this complex can be annihilated by finite covers of the corresponding schemes. By Lemma \ref{extendcover}, we may refine these finite covers by one that comes from all of $S$. In other words, we can find a finite surjective cover $S' \to S$ such that $\alpha'|_{S'} = 0$. After replacing $S$ with $S'$, the Cech spectral sequence then implies that $\alpha$ comes from some $H^{n-q'}(\calU,H^{q'}(G))$ with $q' < q$. Proceeding in this manner, we can reduce the second index $q$ all the way down to $0$, i.e., assume that the class $\alpha$ lies in the image of the map
\[ H^n(\calU,G) \to H^n(S,G). \]
Now we are reduced to the situation in Zariski cohomology that was tackled in Proposition \ref{killcohzarffgs}.
\end{proof}

\begin{remark}
The proof given above for Theorem \ref{killcohffgs} used the intermediary of abelian schemes to make the connection between fppf cohomology and \'etale cohomology with coefficients in a finite flat group commutative scheme $G$ (see Proposition \ref{ffgscohshring}). When the coefficient group scheme $G$ is smooth (or equivalently \'etale), this reduction follows directly from Grothendieck's theorem, thereby giving an elementary proof of Theorem \ref{killcohffgs} in this case. However, we do not know an elementary proof (one that avoids projective geometry) that works in general.
\end{remark}

\begin{remark}
If $G$ is a finite flat group scheme over $S$ which is not necessarily abelian, the $H^1$ part of Theorem \ref{killcohffgs} remains valid since one can trivialise a $G$-torsor $\pi:T \to S$ using the finite flat morphism $\pi$.
\end{remark}

\section{The theorem for abelian schemes}
\label{sec:abschresult}

Our goal in this section is to prove the first half of Theorem \ref{killcohabs}. The arguments here essentially mirror those for finite flat commutative group schemes presented in \S \ref{sec:ffgsresult}. The key difference is that annihilating Zariski cohomology requires more complicated constructions when the coefficients are abelian schemes. We handle this by proving a generalisation of Weil's extension lemma (see Proposition \ref{abschpurity}). This generalisation requires strong regularity assumptions on $S$ and is one of the two places in our proof of Theorem \ref{killcohabs} that we need {\em proper} covers instead of finite ones; the other is the case of $H^1$.

We begin by recording an elementary criterion for a map to an abelian variety to be constant. 

\begin{lemma}
\label{mapcurveabvar}
Let $A$ be an abelian variety over an algebraically closed field $k$, and let $C$ be a reduced variety over $k$. Fix an integer $\ell$ invertible on $k$. A map $g:C \to A$ is constant if and only if it induces the $0$ map $H^1_\et(A,\Q_\ell) \to H^1_\et(C,\Q_\ell)$.
\end{lemma}
\begin{proof}
It suffices to show that a map like $g$ that induces the $0$ map on $H^1$ is trivial. As any $k$-variety is covered by curves, it suffices to show that the map $g$ is constant on all curves in $C$. Thus, we reduce to the case that $C$ is a curve. We may also clearly assume that $C$ is normal, i.e., smooth. Let $\overline{C}$ denote the canonical smooth projective model of $C$. Since $A$ is proper, the map $g$ factors through a map $\overline{g}:\overline{C} \to A$. Since $C$ and $\overline{C}$ are normal, the map $\pi_1(C) \to \pi_1(\overline{C})$ is surjective. Hence, the map $H^1_\et(\overline{C},\Q_\ell) \to H^1_\et(C,\Q_\ell)$ is injective. Thus, to answer the question, we may assume that $C = \overline{C}$ is a smooth projective curve. 

Let $A \hookrightarrow \P^n$ be a closed immersion corresponding to a very ample line bundle $\calL$. The map $g:C \to A$ will be constant if we can show that $g^*\calL$ is not ample, i.e., has degree $0$. As the $\ell$-adic cohomology of an abelian variety is generated in degree $1$ (see \cite[\S 12]{MilneAbVar}), the hypothesis on $H^1$ implies that the map $H^2_\et(A,\Q_\ell) \to H^2(C,\Q_\ell)$ is also $0$. In particular, $g^*(c_1\calL) = 0$, where $c_1(\calL) \in H^2(A,\Q_\ell(1)) \simeq H^2(A,\Q_\ell)$ is the first Chern class of the line bundle $\calL$. Since applying $g^*$ commutes with taking the first Chern class, it follows that $c_1(g^*\calL) = 0$, hence $g^*\calL$ has degree $0$ as desired.
\end{proof}

We now prove the promised extension theorem for maps into abelian schemes.

\begin{proposition}
\label{abschpurity}
Let $S$ be a regular connected excellent noetherian scheme, and let $f:A \to S$ be an abelian scheme. For any non-empty open $U \subset S$, the restriction map $A(S) \to A(U)$ is bijective.
\end{proposition}
\begin{proof}
Let $j:U \to S$ denote the open immersion defined by $U$. The bijectivity of $A(S) \to A(U)$ will follow by taking global sections if we can show that the natural map of presheaves $a:A \to j_*(A|_U)$ is an isomorphism on the small Zariski site of $S$. As both the source and the target of $a$ are actually sheaves for the \'etale topology on $S$, we may localise to assume that $S$ is the spectrum of a strictly henselian local ring $R$. In this setting, we will show that $A(S) \to A(U)$ is bijective using $\ell$-adic cohomology.

The injectivity of $A(S) \to A(U)$ follows from the density of $U \subset S$ and the separatedness of $A \to S$. To show surjectivity, by the valuative criterion of properness, we may assume that the complement $S \setminus U$ has codimension at least $2$ in $S$. Let $s:U \to A$ be a section of $A$ over $U$. By taking the normalised scheme-theoretic closure of $s(U) \subset A$, we obtain a proper birational map $p:S' \to S$ that is an isomorphism over $U$, and an $S$-map $i:S' \to A$ extending $s$ over $U$. The desired surjectivity then reduces to showing that $i$ is constant on the fibres of $p$. Since $p_*\calO_{S'} = \calO_S$, the rigidity lemma (see \cite[Proposition 6.1]{MumfordGIT}) shows that it suffices to show that $i$ collapses the reduced special fibre $S'_s$, where $s \in S$ is the closed point. By Lemma \ref{mapcurveabvar}, it is enough to check that the induced map $H^1(A_s,\Q_\ell) \to H^1(S'_s,\Q_\ell)$ is trivial for some integer $\ell$ invertible on $S$. Note that we have the following commutative diagram:
\[ \xymatrix{ H^1(A,\Q_\ell) \ar[r] \ar[d] & H^1(A_s,\Q_\ell) \ar[d] \\
   H^1(S',\Q_\ell) \ar[r] & H^1(S'_s,\Q_\ell). } \]
The horizontal maps are isomorphisms by the proper base change theorem in \'etale cohomology (see \cite[Arcata IV-1, Th\'eor\`eme 1.2]{SGA4.5}) as $S$ is a strictly henselian local scheme. Hence, it suffices to show that $H^1(A,\Q_\ell) \to H^1(S',\Q_\ell)$ is $0$. Since $H^1(S',\Q_\ell) = \Hom_{\mathrm{conts}}(\pi_1(S'),\Q_\ell)$, it suffices to check that $\pi_1(S') = 0$. As $S'$ is normal, we know that $\pi_1(U) \to \pi_1(S')$ is surjective. Moreover, by Zariski-Nagata purity (see \cite[Expos\'e X, Th\'eor\`eme 3.4]{SGA2}), we know that  $\pi_1(U) \simeq \pi_1(S)$ since $S \setminus U$ has codimension $\geq 2$ in $S$. Since $S$ is strictly henselian, we have $\pi_1(S) = 0$ and hence $\pi_1(S') = 0$ as desired.
\end{proof}

\begin{remark}
\label{rem:topabspurity}
The main idea for the proof of Proposition \ref{abschpurity} comes from obstruction theory in topology. Consider the universal family $\pi:\calU_g \to \calA_g$ of abelian varieties over the stack $\calA_g$ of abelian varieties. Proposition \ref{abschpurity} can be rephrased as asking if every map $S \to \calA_g$ with a specified lift $U \to \calU_g$ over a dense open $U \subset S$ admits an extension $S \to \calU_g$ provided $S$ is smooth. Since the stack $\calU_g$ is a classifying space for its fundamental group (since the same is true for $\calA_g$ and the fibres of $\pi$), the answer at the level of homotopy types would be yes if and only if $\pi_1(U) \to \pi_1(\calU_g)$ factors through $\pi_1(U) \to \pi_1(S)$. This is essentially what is verified above using purity; Lemma \ref{mapcurveabvar} allows us to go from this homotopy-theoretic conclusion to a geometric one.
\end{remark}

\begin{remark}
Proposition \ref{abschpurity} can be considered a generalisation of Weil's extension lemma when applied to abelian varieties. Recall that this lemma says that the domain of definition of rational maps from a smooth variety to a group variety has pure codimension $1$. In case the target is proper, i.e., an abelian variety $A$, this reduces to the statement that $A(X) \simeq A(U)$ for any smooth variety $X$, and dense open $U \hookrightarrow X$.
\end{remark}

\begin{remark}
\label{abspurityratsing}
Our proof of Proposition \ref{abschpurity} is topological as explained in Remark \ref{rem:topabspurity}. As pointed out to us by J\'anos Koll\'ar after the present work was completed, one can also give a more geometric proof of Proposition \ref{abschpurity} as follows: a theorem of Abhyankar (see  \cite[\S VI.1, Theorem 1.2]{KollarRCbook}) implies that for any proper modification $p:S' \to S$ with $S$ noetherian regular excellent, the positive dimensional fibres of $p$ contain non-constant rational curves. Applying this theorem to the closure $S'$ of the graph of a rational map defined by a section $U \to A$ over an open $U \subset S$ gives our desired claim as abelian varieties do not contain rational curves. We prefer the cohomological approach as a slight variation on it (using cohomology of the structure sheaf $\calO_S$ instead of the constant sheaf in the proof of Theorem \ref{abschpurity} and Lemma \ref{mapcurveabvar})  shows that Proposition \ref{abschpurity} remains valid in characteristic $0$ if $S$ has rational singularities. This also suggests a question to which we do not know the answer: if $S$ is a scheme in positive characteristic satisfying some definition of rational singularities (such $F$-rationality), does Proposition \ref{abschpurity} hold for $S$?
\end{remark}

\begin{example}
\label{abspurityfails}
We give an example to show that the regularity condition on $S$ cannot be weakened too much in Proposition \ref{abschpurity}. Let $(E,e) \subset \P^2$ be an elliptic curve, and let $S$ be the affine cone on $E$ with origin $s$. Note that $S$ is a hypersurface singularity of dimension $2$ with $0$ dimensional singular locus. In particular, it is normal. Let $A = S \times E$ denote the constant abelian scheme on $E$ over $S$. Then $U = S \setminus \{s\}$ can be identified with the total space of the $\G_m$-torsor $\calO(-1)|_E - 0(E)$ over $E$. Thus, there exists a non-constant section of $A(U)$. On the other hand, all sections $S \to A$ are constant. Indeed, every point in $S$ lies on an $\A^1$ containing $s$. As all maps $\A^1 \to E$ are constant, the claim follows. Thus, we obtain an example of a normal hypersurface singularity $S$ and an abelian scheme $A \to S$ such that the conclusion of Proposition \ref{abschpurity} fails for $S$. Of course, $S$ is not a rational singularity, a fact supported by Remark \ref{abspurityratsing}.
\end{example}

Next, we point out how to use Proposition \ref{abschpurity} to prove the version of Theorem \ref{killcohabs} involving Zariski cohomology under strong regularity assumptions on the base scheme $S$; the proof is trivial.

\begin{corollary}
\label{killzarcohabs}
Let $S$ be a regular excellent noetherian scheme, and let $f:A \to S$ be an abelian scheme. Then $H^n_\zar(S,A) = 0$ for $n > 0$.
\end{corollary}
\begin{proof}
By Proposition \ref{abschpurity}, we know that $A$ restricts to a constant sheaf on the small Zariski site of each connected component of $S$. By the vanishing of the cohomology of a constant sheaf on an irreducible topological space, the claim follows.
\end{proof}

We are now in a position to complete the proof of Theorem \ref{killcohabs}.

\begin{proof}[Proof of Theorem \ref{killcohabs}]
Let $S$ be a noetherian excellent scheme, and let $A \to S$ be an abelian scheme. We will show that that cohomology classes in $H^n(S,A)$ are killed by proper surjective maps by induction on $n$ provided $n > 0$. We may assume that $S$ is integral. 

For $n = 1$, classes in $H^1(S,A)$ are represented by \'etale $A$-torsors $T$ over $S$. As $T$ is an fppf $S$-scheme, there exists a quasi-finite dominant morphism $U \to S$ such that $T(U)$ is non-empty. By picking an $S$-map $U \to T$ and taking the closure of the image, we obtain a proper surjective cover $S' \to S$ such that $T(S')$ is not empty. This implies that the cohomology class associated to $T$ dies on passage to $S'$, proving the claim.

We next proceed exactly as in the proof of Theorem \ref{killcohffgs} to reduce down to the case of Zariski cohomology. The only difference is that the references to Proposition \ref{ffgscohshring} are replaced by references to Grothendieck's theorem (see \cite[Th\'eor\`eme 11.7]{GrothDixExpIII}) which, in particular, implies that cohomology classes in $H^n(S,A)$ trivialise over an \'etale cover; we omit the details. 

To show the claim for Zariski cohomology, assume first that $S$ is of finite type over $\Z$. In this case, thanks to de Jong's theorems from \cite{dJAlt2}, we can find a proper surjective cover of $S$ with regular total space. Passing to this cover and applying Corollary \ref{killzarcohabs} then solves the problem. In the case that $S$ is no longer of finite type over $\Z$, we reduce to the finite type case using approximation. Indeed, the data $(S,A,\alpha)$ comprising of the base scheme $S$, the abelian scheme $A \to S$, and a Zariski cohomology class $\alpha \in H^n_\zar(S,A)$ can be approximated by similar data with all schemes involved of finite type over $\Z$. Given such an approximating triple $(S',A',\alpha')$ with $S'$ of finite type over $\Z$, we can find a proper surjective map $S'' \to S'$ killing $\alpha'$ by the earlier argument. By functoriality, the pullback $S'' \times_{S'} S \to S$ is a proper surjective cover of $S$ killing $\alpha$.
\end{proof}

\begin{remark}
Theorem \ref{killcohabs} admits a topological reformulation as follows. Given a noetherian scheme $S$ and an abelian scheme $G$ over $S$, let $S_\prop$, $S_\fppf$ and $S_{\prop,\fppf}$ denote the (big) topoi associated to the category of schemes over $S$ equipped with the topology generated respectively by proper surjective maps, fppf maps, and both proper surjective and fppf maps. There are natural forgetful maps $a:S_{\prop,\fppf} \to S_\prop$ and $b:S_{\prop,\fppf} \to S_\fppf$ of topoi. Given a finite flat commutative group scheme $G \to S$, let $G$ also denote the sheaf defined by $G$ in each of the above topologies. Theorem \ref{killcohabs} can be reformulated as saying that the sheaves $\R^i a_* G$ vanish for $i > 0$. Since schemes are sheaves for the fppf topology, one can easily show that $a_* a^* G = G$. Thus, Theorem \ref{killcohabs} can be reformulated saying that $G \simeq \R a_* G$.  Note a consequence: since cohomology on a site is computed using hypercovers by Verdier's theorem (see \cite[Theorem 8.16]{ArtinMazur}), the previous identification gives that for a class $\alpha \in H^n(S_\fppf,G)$ with $n > 0$,  there exists a proper hypercovering $f_\bullet:T_\bullet \to S$ and a map of simplicial schemes $\phi:T_\bullet \to K(G,n)$ representing $b^* \alpha$.  If $G$ is instead a finite flat group scheme, then the same remarks apply for Theorem \ref{killcohffgs}, except that we replace proper maps by finite ones.
\end{remark}

\section{An application: big Cohen-Macaulay algebras in positive characteristic}
\label{sec:altproof}

Let $(R,\fm)$ be an excellent noetherian local domain containing $\F_p$. A fundamental theorem of Hochster-Huneke (see \cite{HHBigCM}) asserts that the absolute integral closure $R^+$ (the integral closure of $R$ in a fixed algebraic closure of its fraction field) is a Cohen-Macaulay algebra. This result and the ideas informing it form the bedrock of tight closure theory and huge swathes of positive characteristic commutative algebra. 

Our goal in this section is to give a new proof of the Hochster-Huneke theorem using Theorem \ref{killcohffgs}. We hasten to remark that there already exist alternative proofs in the literature, all cocycle-theoretic or equational at the core. The approach adopted here follows closely the relatively-recent approach from \cite{GLBigCM}, the essential new feature being the use of cohomology-annihilation result proven in Theorem \ref{killcohffgs} in place of explicit cocycle manipulations.

We begin by recording a coherent cohomology-annihilation result one can deduce from Theorem \ref{killcohffgs}; this can be considered as the analogue  of the ``equational lemma'' of \cite{HHBigCM}; see also \cite[Lemma 2.2]{GLBigCM}.

\begin{proposition}
\label{killlocalcoh}
Let $(R,\fm)$ be a noetherian excellent local $\F_p$-algebra, and let $M \subset H^i_\fm(R)$ be a Frobenius stable finite length $R$-submodule for some $i > 0$. Then there exists a module-finite extension $f:R \to S$ such that $f^*(M) = 0$ where $f^*:H^i_\fm(R) \to H^i_\fm(S)$ is the induced map.
\end{proposition}
\begin{proof}
After normalising $R$, we may assume that $i > 1$. With $U = \Spec(R) - \{\fm\}$, we have a Frobenius equivariant identification
\[ H^{i-1}(U,\calO) \simeq H^i_\fm(R) \]
which allows us to view $M$ as a submodule of $H^{i-1}(U,\calO)$. The Frobenius action endows $H^{i-1}(U,\calO)$ with the structure of a $R\{X^p\}$-module, where $R\{X^p\}$ is the non-commutative polynomial ring over $R$ with one generator $X^p$ satisfying the commutation relation $X^p r = r^p X^p$ for $r \in R$. The finite length assumption implies that for each $m \in M$, there exists some monic additive polynomial $g(X^p) \in R\{X^p\}$ such that $g(m) = 0$. As $g$ is additive and monic, we have a short exact sequence 
\[ 0 \to \ker(g) \to \calO \stackrel{g}{\to} \calO \to 0\]
of abelian sheaves on $\Spec(R)_\fppf$. Moreover, the monicity of $g$ also shows that the sheaf $\ker(g)$ is representable by a finite flat commutative group scheme over $\Spec(R)$. As $g(m) = 0$,  we see that $m$ comes from a cohomology class $m' \in H^{i-1}(U,\ker(g))$. Since $i-1 > 0$, Theorem \ref{killcohffgs} shows that there exists a finite surjective map $\pi:V \to U$ such that $\pi^* m' = 0$. Setting $S$ to be the (global sections of the) normalisation of $R$ in $V$ is then seen to solve the problem.
\end{proof}

Using Proposition \ref{killlocalcoh}, we can give a proof that $R^+$ is Cohen-Macaulay. The argument given below is based entirely on \cite[Theorem 2.1]{GLBigCM} and simply recorded here for convenience.

\begin{theorem}
Let $(R,\fm)$ be a noetherian excellent local $\F_p$-domain that admits a dualising complex, and let $R^+$ be its absolute integral closure. Then $R^+$ is Cohen-Macaulay.
\end{theorem}
\begin{proof}
For notational simplicity, we restrict ourselves to the case that $R$ is complete. Let $d = \dim(R)$, and for a prime $\fp \in \Spec(R)$, we let $d_{\fp} = \dim(R_{\fp})$. For any integer $i$, the Matlis dual $D(H^i_\fm(R))$ of $H^i_\fm(R)$ is a finitely generated $R$-module whose stalks can be described by the formula
\[ D(H^i_\fm(R))_\fp = D(H^{i - (d - d_\fp)}_\fp(R_\fp)). \]
The important point here is that the cohomological degree drops by the codimension on localising (and dualising); see \cite[Expos\'e VIII, Th\'eor\`eme 2.1]{SGA2} for a nice application of this observation.

Let us now show by induction on $d$ that there exists a module-finite extension $R \to S$ which kills all local cohomology outside degree $d$; clearly this suffices to prove the theorem. The case $d = 0$ being vacuous, we assume $d > 0$ and pick a non-negative integer $i < d$. The Matlis dual $D(H^i_\fm(R))$ is a finitely generated $R$-module and hence has a finite set $\fp_1,\dots,\fp_n$ of non-maximal associated primes. For each such prime $\fp_j$, induction constructs a module-finite extension $R_{\fp_j} \to S_j$ that kills $H^{i - (d - d_{\fp_j})}_{\fp_j}(R_{\fp_j})$; note that $i - (d - d_{\fp_j}) < d_{\fp_j}$ since $i < d$. Setting $S$ to be the normalisation of $R$ in a compositum of all the $S_j$'s then shows that the map $R \to S$ induces a map $f_*:D(H^i_\fm(S)) \to D(H^i_\fm(R))$ whose image is not supported at any $\fp_j$. Since the only other possible associated prime of $D(H^i_\fm(R))$ is $\fm$, the support of $\im(f_*)$ is a finite length submodule of $D(H^i_\fm(R))$. By duality, the image of $f^*:H^i_\fm(R) \to H^i_\fm(S)$ is a finite length $R$-submodule $M$ of $H^i_\fm(S)$ which is moreover Frobenius stable. Hence, the $S$-submodule of $H^i_\fm(S) \simeq H^i_{\fm S}(S)$  generated by $M$ is also Frobenius stable with finite length. Since $\fm S$ is a finite colength ideal in $S$, Proposition \ref{killlocalcoh} gives a module-finite extension $S \to T$ killing $M$; the composite $R \to T$ then kills $H^i_\fm(R)$.
\end{proof}

\section{An example of a torsor not killed by finite covers}
\label{abspurityneedh}

Theorem \ref{killcohabs} lets one construct proper covers annihilating cohomology with coefficients in an abelian scheme. Our goal in this section is to construct an example indicating why ``proper'' cannot be replaced by ``finite'' in preceding statement; the key geometric idea in this construction (that of passing to finite \'etale covers of {\em Zariski} local rings) comes from \cite[Example 3.2, Chapter XIII]{RaynaudAmpleLBLNM}. As a bonus, we get an example illustrating the necessity of strong regularity assumptions on the base scheme in Theorem \ref{killzarcohabs}.

\subsection{Construction}
\label{absex:cons}
Fix an algebraically closed field $k$ of characteristic $0$, and an elliptic curve $(E,0)$ over $k$. We will construct a scheme $X$ essentially of finite type over $k$ satisfying the following:
\begin{enumerate}
\item $X$ is a semilocal, normal, connected, $2$-dimensional affine scheme with two closed points $x$ and $y$. Let $U = X - \{x,y\}$ be the twice-punctured spectrum; let $X_x = \Spec(\calO_{X,x})$ and $X_y = \Spec(\calO_{X,y})$ be the corresponding local rings; let $U_x = U \times_X X_x$ and $U_y = U \times_X X_y$ denote the corresponding punctured spectra; and let $\widehat{U_x} = U_x \times_{X_x} \widehat{X_x}$ and $\widehat{U_y} = U_y \times_{X_y} \widehat{X_y}$ denote the punctured spectra of the corresponding completions.
\item All maps $X_x \to E$ induce the trivial map $\pi_1(\widehat{U_x}) \to \pi_1(E)$. 
\item All maps $X_y \to E$ induce the trivial map $\pi_1(\widehat{U_y}) \to \pi_1(E)$.
\item There exists a map $f:U \to E$ inducing surjective maps $\pi_1(\widehat{U_x}) \to \pi_1(E)$ and $\pi_1(\widehat{U_y}) \to \pi_1(E)$ simultaneously.
\end{enumerate}

We first explain the idea of the construction informally. The cone $S$ considered in Example \ref{abspurityfails} had the property that sections of $E$ on an open subscheme do not extend to the entire scheme. Glueing two such cones away from the cone point gives an $E$-torsor of infinite order on a normal scheme that does not die on passage to finite covers. The base scheme, however, is not separated. To achieve separatedness, instead of glueing naively, we look at a finite \'etale quadratic cover of $S$. The resulting scheme bares enough formal similarities with the non-separated example (namely, exactly two closed points, each of which looks like the cone $S$) to make this construction work. The details follow; we advise the reader willing to take the existence of $X$ on faith to proceed to \S \ref{absex:verify}.

Let $\overline{S}$ the affine cone on $E$ (considered in Example \ref{abspurityfails} with different notation), and let $(S,s)$ denote its local scheme at the origin. Note that $S$ is a normal, Gorenstein, local $2$-dimensional scheme essentially of finite type over $k$. By Noether normalisation, we can pick a finite map $a:S \to \Spec(A)$ where $A$ is the local ring at the origin of $\A^2_k$. Since $S$ is Cohen-Macaulay, the map $f$ is finite flat. In fact, we can even arrange for $f$ to be totally ramified at the origin: if $E$ is represented by the homogeneous form $y^2z = x^3 + Axz^2 + Bz^3$, then we simply choose $a$ to be the map given by the functions $y$ and $z$. Let $b:\Spec(B) \to \Spec(A)$ be a finite \'etale cover of degree $2$ with $B$ connected. We define $X$ to be the fibre product via 
\[ \xymatrix{ X := S \times_{\Spec(A)} \Spec(B) \ar[d] \ar[r] & \Spec(B) \ar[d]^b \\ S \ar[r]^a & \Spec(A). } \]
The scheme $X$ is connected since $b$ is \'etale at the origin while $a$ is totally ramified. Moreover, being finite \'etale over $S$ forces $X$ to be a semilocal, normal, connected, $2$-dimensional affine scheme with two closed points $x$ and $y$. Note that since $X \to S$ is finite \'etale of degree $2$, it is necessarily Galois. Let $X_x$, $X_y$, $U_x$, etc. be as above, and let's verify the desired properties.

First, we verify properties $(2)$ and $(3)$. The map $\widehat{U_x} \to E$ induced by a map $X_x \to E$ factors through the induced map $\widehat{X_x} \to E$ by definition. Since $\widehat{X_x}$ is a complete noetherian local scheme with algebraically closed residue field, its fundamental group vanishes, and hence the desired claim follows for $X_x$. We argue exactly the same way for $X_y$.

For property $(4)$, we first explain how to construct $f$. Consider the punctured spectrum $S - \{s\}$. As explained in Example \ref{abspurityfails}, the scheme $S - \{s\}$ can be realised as the complement of the $0$-section in the Zariski localisation along the $0$-section of the total space of the line bundle $\calO(-1)|_E \to E$. In particular, there exists a natural map $f_0:S - \{s\} \to E$. Let $f$ denote the composition $U \to S - \{s\} \to E$. Note that $f$ is invariant under the Galois group of $X \to S$.

We will now verify that $f$ has the desired properties. Note that since $S - \{s\}$ is ind-open in $\calO(-1)|_E$ and both schemes are normal, the induced map $\pi_1(S - \{s\}) \to \pi_1(\calO(-1)|_E)$ is surjective. Since we are working in characteristic $0$, by homotopy invariance of the fundamental group for normal schemes, we can identify $\pi_1(\calO(-1)|_E) \simeq \pi_1(E)$ via the natural projection. In particular, $f_0$ induces a surjective map $\pi_1(S - \{s\}) \to \pi_1(E)$. Moreover, the same calcuations also work after completion. Hence, the induced map $\widehat{f_0}:\widehat{S - \{s\}} \to E$ also induces a surjective map on fundamental groups. To pass to $X$, note that $X \to S$ is finite \'etale with no residue extension at the closed points. Hence, the induced map $\widehat{X_x} \to \widehat{S}$ is an isomorphism, and similarly for $y$. This allows us to identify $\widehat{U_x}$ with $\widehat{S - \{s\}}$ via the natural map, and similarly for $y$. The desired surjectivity now follows from what we already checked for $S$.

\subsection{Verification}
\label{absex:verify}

Let $X$ be the scheme constructed in \S \ref{absex:cons}. Let $V_x = X - \{y\}$ and $V_y = X - \{y\}$ denote the open subschemes of $X$ defined as the complement of each of the two closed points. Since $x$ and $y$ are the only two closed points of $X$, the pair $\{V_x,V_y\}$ defines a Zariski open cover of $X$ with intersection $U = V_x \cap V_y$. Consider the associated Mayer-Vietoris sequence
\begin{equation}\label{absex:exseq} \cdots H^0(V_x,E) \oplus H^0(V_y,E) \stackrel{\beta}{\to} H^0(U,E) \stackrel{\delta}{\to} H^1(X,E) \to \cdots \end{equation}

By assumption, we may pick a map $f:U \to E$ inducing surjective maps $\pi_1(\widehat{U_x}) \to \pi_1(E)$ and $\pi_1(\widehat{U_y}) \to \pi_1(E)$. Viewing $f$ as an element of $H^0(U,E)$, we define $\alpha = \delta(f)$. We claim:

\begin{claim}
\label{absex:claim}
The class $\alpha$ is not torsion. Moreover, for every finite cover $f:Y \to X$, the pullback $f^*\alpha$ is also not torsion. In particular, $\alpha$ does not die in a finite cover.
\end{claim}

\begin{proof}
Assuming that $\alpha$ is not torsion, the existence of norms (see Proposition \ref{absex:norms}) implies the rest. Thus, it suffices to verify that $\alpha = \delta(f)$ is not torsion. By exact sequence (\ref{absex:exseq}), it suffices to show that $n \cdot f \notin \im(\beta)$ for an integer $n \neq 0$. The class $n \cdot f$ is represented by the map $[n] \circ f:U \to E$ where $[n]:E \to E$ is the multiplication by $n$ map on $E$. The assumption on $f$ implies that $[n] \circ f$ induces a map $\pi_1(\widehat{U_x}) \to \pi_1(E)$ whose image has finite index. On the other hand, our assumptions on $X$ also imply that any map lying in the image of $\beta$ induces trivial map $\pi_1(\widehat{U_x}) \to \pi_1(E)$. Since $\pi_1(E)$ is torsion free and non-zero, the claim follows.
\end{proof}

Lastly, we explain why the functors $H^i(-,A)$ admit ``norm maps'' when $A$ is a group algebraic space, and why this forces cohomology classes killed by finite covers to be torsion. These maps were used above.

\begin{proposition}
\label{absex:norms}
Let $A$ be any abelian fppf sheaf on the category of all $k$-schemes with $k$ a field of characteristic $0$. Assume that $A$ is represented by an algebraic space. For every finite surjective morphism $f:T \to S$ of integral normal schemes, there exist norm maps $f_*:A(T) \to A(S)$ satisfying the conditions listed in \cite[Definition 4.1]{SuslinVoevodskySingHom}. Moreover, the kernel of $H^i(S,A) \to H^i(T,A)$ is necessarily torsion.
\end{proposition}

Proposition \ref{absex:norms} is well-known, but we were unable to find a satisfactory reference. Hence, we include the sketch of a proof.

\begin{proof}[Sketch of proof]
We first explain the construction of norms. Assume that $T \to S$ induces a Galois extension of function fields with group $G$ with cardinality $n$. By normality of $S$, we identify $T/G \simeq S$, where the quotient $T/G$ is computed in the category of algebraic spaces (or schemes). Given a $T$-point $a \in A(T)$, we obtain a natural map $T \to \Map(G,A) \simeq A^n$ given by $t \mapsto (g \mapsto a(g(t)))$. The group $S_n = S_{\#G}$ acts on $\Map(G,A)$, and the preceding map $T \to \Map(G,A)$ is equivariant for the natural embedding $G \to S_n$ given by left translation. Taking quotients as algebraic spaces, we arrive at a map $b:S \simeq T/G \to A^n/S_n = \Sym^n(A)$. The $n$-fold multiplication map $A^n \to A$ is an $S_n$-equivariant map to an algebraic space. Hence, it factors as $A^n \to \Sym^n(A) \to A$. Composing the second map with $b$, we obtain a map $S \to A$ that we declare to be the norm $f_*(a)$ of $a$. In the case $f:T \to S$ does not induce a Galois extension of function fields, one works in a Galois closure and then descends; we omit the details as they do not mattter in the sequel. One can check that this constructions verifies the conditions required on the norm map in \cite[Definition 4.1]{SuslinVoevodskySingHom}. 

For the last claim, note that a formal consequence of having norms is that for any finite surjective morphism $f:T \to S$ as above, there exists a pushforward $H^i(f_*):H^i(T,A) \to H^i(S,A)$ such that $H^i(f_*) \circ H^i(f^*) = n$. In particular, a non-torsion class in $H^i(S,A)$ will not die on passage to finite covers if $S$ is normal. 
\end{proof}

\bibliography{mybib}

\begin{thebibliography}{MFK94}

\bibitem[AM69]{ArtinMazur}
M.~Artin and B.~Mazur.
\newblock {\em Etale homotopy}.
\newblock Lecture Notes in Mathematics, No. 100. Springer-Verlag, Berlin, 1969.

\bibitem[Art74]{ArtinVersalDefs}
M.~Artin.
\newblock Versal deformations and algebraic stacks.
\newblock {\em Invent. Math.}, 27:165--189, 1974.

\bibitem[BBM82]{BBMII}
Pierre Berthelot, Lawrence Breen, and William Messing.
\newblock {\em Th\'eorie de {D}ieudonn\'e cristalline. {II}}, volume 930 of
  {\em Lecture Notes in Mathematics}.
\newblock Springer-Verlag, Berlin, 1982.

\bibitem[Con07]{ConradNagata}
Brian Conrad.
\newblock Deligne's notes on {N}agata compactifications.
\newblock {\em J. Ramanujan Math. Soc.}, 22(3):205--257, 2007.

\bibitem[Del77]{SGA4.5}
Pierre Deligne.
\newblock {\em Cohomologie \'etale}.
\newblock Springer-Verlag, Berlin, 1977.
\newblock S\'eminaire de G\'eom\'etrie Alg\'ebrique du Bois-Marie SGA 4${1\over
  2}$, Avec la collaboration de J. F. Boutot, A. Grothendieck, L. Illusie et J.
  L. Verdier, Lecture Notes in Mathematics, Vol. 569.

\bibitem[dJ97]{dJAlt2}
A.~J. de~Jong.
\newblock Families of curves and alterations.
\newblock {\em Ann. Inst. Fourier (Grenoble)}, 47(2):599--621, 1997.

\bibitem[FC90]{FaltingsChai}
Gerd Faltings and Ching-Li Chai.
\newblock {\em Degeneration of abelian varieties}, volume~22 of {\em Ergebnisse
  der Mathematik und ihrer Grenzgebiete (3) [Results in Mathematics and Related
  Areas (3)]}.
\newblock Springer-Verlag, Berlin, 1990.
\newblock With an appendix by David Mumford.

\bibitem[Fon94]{Fontainepadicperiod}
Jean-Marc Fontaine.
\newblock Le corps des p\'eriodes {$p$}-adiques.
\newblock {\em Ast\'erisque}, (223):59--111, 1994.
\newblock With an appendix by Pierre Colmez, P{\'e}riodes $p$-adiques
  (Bures-sur-Yvette, 1988).

\bibitem[Gro66]{EGA4_3}
A.~Grothendieck.
\newblock \'{E}l\'ements de g\'eom\'etrie alg\'ebrique. {IV}. \'{E}tude locale
  des sch\'emas et des morphismes de sch\'emas. {III}.
\newblock {\em Inst. Hautes \'Etudes Sci. Publ. Math.}, (28):255, 1966.

\bibitem[Gro67]{EGA4_4}
A.~Grothendieck.
\newblock \'{E}l\'ements de g\'eom\'etrie alg\'ebrique. {IV}. \'{E}tude locale
  des sch\'emas et des morphismes de sch\'emas {IV}.
\newblock {\em Inst. Hautes \'Etudes Sci. Publ. Math.}, (32):361, 1967.

\bibitem[Gro68a]{SGA2}
Alexander Grothendieck.
\newblock {\em Cohomologie locale des faisceaux coh\'erents et th\'eor\`emes de
  {L}efschetz locaux et globaux {$(SGA$} {$2)$}}.
\newblock North-Holland Publishing Co., Amsterdam, 1968.
\newblock Augment{\'e} d'un expos{\'e} par Mich{\`e}le Raynaud, S{\'e}minaire
  de G{\'e}om{\'e}trie Alg{\'e}brique du Bois-Marie, 1962, Advanced Studies in
  Pure Mathematics, Vol. 2.

\bibitem[Gro68b]{GrothDixExpIII}
Alexander Grothendieck.
\newblock Le groupe de {B}rauer. {III}. {E}xemples et compl\'ements.
\newblock In {\em Dix {E}xpos\'es sur la {C}ohomologie des {S}ch\'emas}, pages
  88--188. North-Holland, Amsterdam, 1968.

\bibitem[Gro03]{SGA1}
A.~Grothendieck.
\newblock {\em Rev\^etements \'etales et groupe fondamental ({SGA} 1)}.
\newblock Documents Math\'ematiques (Paris) [Mathematical Documents (Paris)],
  3. Soci\'et\'e Math\'ematique de France, Paris, 2003.

\bibitem[HH92]{HHBigCM}
Melvin Hochster and Craig Huneke.
\newblock Infinite integral extensions and big {C}ohen-{M}acaulay algebras.
\newblock {\em Ann. of Math. (2)}, 135(1):53--89, 1992.

\bibitem[HL07]{GLBigCM}
Craig Huneke and Gennady Lyubeznik.
\newblock Absolute integral closure in positive characteristic.
\newblock {\em Adv. Math.}, 210(2):498--504, 2007.

\bibitem[Hoc07]{HochsterSurvey}
M.~Hochster.
\newblock Homological conjectures, old and new.
\newblock {\em Illinois J. Math.}, 51(1):151--169 (electronic), 2007.

\bibitem[Hoo82]{HooblerGabber}
Raymond~T. Hoobler.
\newblock When is {${\rm Br}(X)={\rm Br}\sp{\prime} (X)$}?
\newblock In {\em Brauer groups in ring theory and algebraic geometry
  ({W}ilrijk, 1981)}, volume 917 of {\em Lecture Notes in Math.}, pages
  231--244. Springer, Berlin, 1982.

\bibitem[Kol96]{KollarRCbook}
J{\'a}nos Koll{\'a}r.
\newblock {\em Rational curves on algebraic varieties}, volume~32 of {\em
  Ergebnisse der Mathematik und ihrer Grenzgebiete. 3. Folge. A Series of
  Modern Surveys in Mathematics [Results in Mathematics and Related Areas. 3rd
  Series. A Series of Modern Surveys in Mathematics]}.
\newblock Springer-Verlag, Berlin, 1996.

\bibitem[LMB00]{LMBChAlg}
G{\'e}rard Laumon and Laurent Moret-Bailly.
\newblock {\em Champs alg\'ebriques}, volume~39 of {\em Ergebnisse der
  Mathematik und ihrer Grenzgebiete. 3. Folge. A Series of Modern Surveys in
  Mathematics [Results in Mathematics and Related Areas. 3rd Series. A Series
  of Modern Surveys in Mathematics]}.
\newblock Springer-Verlag, Berlin, 2000.

\bibitem[MFK94]{MumfordGIT}
D.~Mumford, J.~Fogarty, and F.~Kirwan.
\newblock {\em Geometric invariant theory}, volume~34 of {\em Ergebnisse der
  Mathematik und ihrer Grenzgebiete (2) [Results in Mathematics and Related
  Areas (2)]}.
\newblock Springer-Verlag, Berlin, third edition, 1994.

\bibitem[Mil]{MilneAbVar}
J.~Milne.
\newblock Abelian {V}arieties.
\newblock Available at
  \url{http://http://www.jmilne.org/math/CourseNotes/av.html}.

\bibitem[Ray70]{RaynaudAmpleLBLNM}
Michel Raynaud.
\newblock {\em Faisceaux amples sur les sch\'emas en groupes et les espaces
  homog\`enes}.
\newblock Lecture Notes in Mathematics, Vol. 119. Springer-Verlag, Berlin,
  1970.

\bibitem[SV96]{SuslinVoevodskySingHom}
Andrei Suslin and Vladimir Voevodsky.
\newblock Singular homology of abstract algebraic varieties.
\newblock {\em Invent. Math.}, 123(1):61--94, 1996.

\end{thebibliography}

\end{document}